\documentclass[a4paper,reqno]{amsart}

\usepackage{amsmath}
\usepackage{amsthm}
\usepackage{amsfonts}
\usepackage{amssymb}

\usepackage{color}

\usepackage{hyperref}

\hypersetup{
    colorlinks=false,
    pdfborder={0 0 0},
    pdfstartview={XYZ null null 1.00}
}

\numberwithin{equation}{section}

  \def\<{\langle}
  \def\>{\rangle}

\theoremstyle{plain}
  \newtheorem{theorem}{Theorem}[section]
  \newtheorem{proposition}[theorem]{Proposition}
  \newtheorem{lemma}[theorem]{Lemma}

\theoremstyle{definition}
  \newtheorem{definition}[theorem]{Definition}

\begin{document}

\title[Krasnosel'skii type formula...]{Krasnosel'skii type formula and translation along trajectories method on the scale of fractional spaces}

\author{Piotr Kokocki}
\address{\noindent BCAM - Basque Center for Applied Mathematics \newline Alameda de Mazarredo 14, 48009, Bilbao, Spain}
\address{\noindent  Faculty of Mathematics and Computer Science \newline Nicolaus Copernicus University \newline  Chopina 12/18, 87-100 Toru\'n, Poland}
\email{pkokocki@mat.umk.pl}
\thanks{The researches supported by the NCN Grant no. NCN  2013/09/B/ST1/01963}



\begin{abstract}
We provide global continuation principle of periodic solutions for the equation $\dot u = - Au  + F(t,u)$, where $A:D(A)\to X$ is a sectorial operator on a Banach space $X$ and $F:[0,+\infty)\times X^\alpha\to X$ is a nonlinear map defined on fractional space $X^\alpha$. The approach that we use in this paper is based upon the theory of topological invariants that applies in the situation when Poincar\'e operator associated with the equation is endowed with some form of compactness.
\end{abstract}

\maketitle

\setcounter{tocdepth}{2}

\section{Introduction}

In this paper we are concerned with $T$-periodic solutions for the differential equation of the following form
\begin{equation}\label{eq-intr}
\dot u(t) = - A u(t) +F (t,u(t)) \quad \text{ on } \ \ [0,+\infty)
\end{equation}
where $A:X\supset D(A)\to X$ is a sectorial operator with compact resolvents on a Banach space $X$ and $F:[0,+\infty)\times X^\alpha\to X$ is a continuous map defined on the fractional space $X^\alpha$ associated with the operator $A$ (see e.g. \cite{MR1721989}, \cite{Henry}, \cite{H-F}, \cite{Pazy} for construction and more details on fractional spaces).

Given $x\in X^\alpha$, let us assume that $u(t;x)$ is a mild solution of the equation \ref{eq-intr} such that $u(0;x) = x$. We look for the $T$-periodic solutions of this equation as the fixed points of the Poincar\'e operator $\Phi_T \colon X^\alpha \to X^\alpha$ given by the formula $\Phi_T(x) := u(T;x)$. Since $A$ has compact resolvent, under usual assumptions on the nonlinearity $F$, one can prove that the operator $\Phi_T$ is completely continuous. Therefore, the natural way to obtain the existence of fixed points for the Poincar\'e operator is to
apply the approach based on topological degree theory that allows us to obtain effective methods to search for the fixed points of the Poincar\'e operator. To be more precise, we will strongly use the homotopy invariance of topological degree to obtain the following results: {\em Krasnosel'skii type degree formula} and {\em averaging principle}, that will lead us to {\em the global continuation principle} of $T$-periodic solutions for the equation \ref{eq-intr}.

The Krasnosel'skii type degree formula states that, if $F$ is time-independent and $U\subset X^\alpha$ is an open bounded set, then $$\mathrm{deg_{LS}}(I - \Phi_t, U) = \mathrm{deg_\alpha}(-A + F, U)$$ for sufficiently small $t>0$. Here $\mathrm{deg_{LS}}$ ia a Leray-Schauder degree on the space $X^\alpha$ and $\mathrm{deg_\alpha}$ is a topological degree for perturbations of sectorial operators (see Appendix for more details). Obtained degree formula is used to derive {\em the averaging principle} that concerns the following family of differential equations
\begin{equation}\label{eewww}
\dot u(t) = - \lambda A u(t) + \lambda F (t,u(t)) \quad \text{ on } \ \ [0,+\infty),
\end{equation}
where $\lambda\in(0,1]$ is a parameter. Let us assume that $\Phi^\lambda_T : X^\alpha\to X^\alpha$ is the Poncar\'e operator associated with the equation \ref{eewww} and let $-A + \widehat{F}:\overline{U}\cap D(A)\to X$ be the map where $$\widehat{F}(x) := \frac{1}{T}\int_0^TF(\tau,x)\,d \tau \qquad\mathrm{for}\quad x\in X^\alpha.$$
{\em The averaging principle} asserts that, for any open bounded set $U\subset X^\alpha$ such that $0\not\in(-A + \widehat{F})(\partial U)$, one choose sufficiently small $\lambda > 0$ such that
\begin{equation}\label{aaddff}
\mathrm{deg_{LS}}(I - \Phi^\lambda_T, U) = \mathrm{deg_\alpha}(-A + \widehat{F}, U).
\end{equation}
The formula \ref{aaddff} allows us to study the periodic solutions for the equation \ref{eewww} only if the positive parameter $\lambda$ is sufficiently close to zero. To improve this for the whole interval $(0,1]$, we prove {\em the global continuation principle}, that provides conditions ensuring us that the topological degree $\mathrm{deg_\alpha}(-A + \widehat{F}, U)$ is non-trivial and, for any $\lambda\in(0,1]$, the Poincar\'e operator $\Phi^\lambda_T$ does not admit fixed points on $\partial U$. This together with the homotopy invariance of topological degree allow us to obtain the existence of $T$-periodic solutions for the equation \ref{eq-intr}. Observe that the advantage of using topological degree is that we obtain not only the existence of periodic solution but also we compute its topological degree. This is an additional information that is crucial in the examination of the stability and multiplicity of periodic solutions. These issues constitute the subject of our subsequent studies.

The Krasnosel'skii type degree formula and averaging principle were considered initially in \cite{Furi-Pera}, \cite{MR819192}, \cite{MR2376467}, \cite{MR0442775}, \cite{MR1368675}, \cite{MR1744938}, \cite{MR2766807} for differential equations on finite dimensional spaces or manifolds. Their generalization on the case of equations on arbitrary Banach spaces were considered in \cite{Cwiszewski-1}, where the right side of the equation is the nonlinear perturbation of the generator of a compact $C_0$ semigroup. See also \cite{Cwiszewski-2} for results in the case when $A$ is a general single valued $m$-accretive operator. In \cite{cw-kok1} the averaging principle were studied in the case when $-A$ generates a $C_0$ semigroup of contractions and $F$ is a condensing map with respect to the Hausdorff measure of noncompactness. Obtained results were used to study the periodic solutions of the first order hyperbolic equations.
The case of weakly damped wave equation were studied in \cite{MR2849823}, while the fourth order beam equation were considered in \cite{MR2870924}.
Generalization of the averaging principle for the equations where the right side is the nonlinear perturbation of the family of generators of $C_0$ semigroups $\{A(t)\}_{t\ge 0}$ is contained in \cite{cwi-kok2}. It is also worth to note that the principle was used in \cite{MR2772422} to prove the existence of positive solutions for nonlinear parabolic equations.

Common feature of the results contained in these papers is the fact that they cover only the case where the domain of $F$ is exactly the space $X^0 = X$. This raises difficulties in applications because the broad class of partial differential equations where derivatives are involved in the nonlinearity do not fit into this setting. The class includes Cahn-Hilliard, Navier-Stokes and parabolic partial differential equations with nonlinearity depending on gradient (for more details see e.g. \cite{DloCho}, \cite{Henry}). In this paper we overcome this problem and prove the new Krasnosel'skii type degree formula and averaging principle for the class of differential equations, where $A$ is the sectorial operator and the domain of the nonlinearity $F$ is any fractional space $X^\alpha$, where $\alpha\in[0,1)$. As an application we provide the global continuation principle of $T$-periodic solutions and use it to prove the criterium for finding $T$-periodic solutions for the second order nonlinear parabolic equations, where the gradient is involved in nonlinearity. We remark that the order of the differential operator is in this case irrelevant
and the results can be similarly applied to equations with higher order differential operators. \\[5pt]

\noindent {\em \bf Notation and terminology.} Let $A:X\supset D(A)\to X$ be a linear operator on a real Banach space $X$ equipped with the norm $\|\cdot\|$. We say that $A$ is sectorial provided there are $\phi\in(0,\pi/2)$,  $M\ge 1$ and $a\in\mathbb{R}$, such that the sector $$S_{a,\phi}:=\{\lambda\in\mathbb{C} \ | \ \phi \le |\mathrm{arg} \, (\lambda - a) \le \pi, \ \lambda\neq a\}$$ is contained in the resolvent set of $A$ and furthermore $$\|(\lambda I - A)^{-1}\| \le M/ |\lambda - a| \qquad\mathrm{for}\quad \lambda\in S_{a,\phi}.$$
It is well-known that $-A$ is an infinitesimal generator of analytic semigroup which, throughout this paper, will be denoted by $\{S_A(t)\}_{t\ge 0}$. The operator $A$ is called positive if $\Re \mu > 0$ for any $\mu\in\sigma(A)$. It can be proved that, if $A$ is positive and sectorial, then given $\alpha \ge 0$ the integral
\begin{equation*}
A^{-\alpha} := \frac{1}{\Gamma(\alpha)}\int_0^\infty t^{\alpha - 1}S_A(t) \,d t.
\end{equation*}
is convergent in the uniform operator topology of the space $\mathcal{L}(X)$. Consequently we can define {\em the fractional space} associated with $A$ as the domain of the inverse operator $X^\alpha:= D(A^\alpha)$. The space $X^\alpha$ endowed with the graph norm $\|x\|_\alpha := \|A^\alpha x\|$ is a Banach space, continuously embedded in $X$. We refer the reader to \cite{Henry}, \cite{H-F}, \cite{Pazy} for more details on sectorial operators and fractional spaces.

\section{Continuity and compactness properties of Cauchy problem}

In this section we consider the following family of differential equations
\begin{equation}\label{row-a-fs}
\dot u(t) = - A u(t) + F (s,t,u(t)), \qquad  t > 0,
\end{equation}
where $s\in[0,1]$ is a parameter, $A: X\supset D(A)\to X$ is a sectorial operator with compact resolvents on a Banach space $X$ and $F:[0,1]\times[0,+\infty)\times X^\alpha\to X$ is a continuous map on fractional space $X^\alpha$, $\alpha\in[0,1)$, such that: \\[5pt]
\noindent\makebox[9mm][l]{$(F1)$}\parbox[t]{118mm}{for every $s\in[0,1]$ and $x\in X^\alpha$ there is an open neighborhood $V\subset X^\alpha$ of $x$ and constant $L > 0$ such that for any $x_1,x_2\in V$ and $t\in [0,+\infty)$ $$\|F(s,t,x_1) - F(s,t,x_2)\|\le L \|x_1 - x_2\|_\alpha;$$}\\
\noindent\makebox[9mm][l]{$(F2)$}\parbox[t]{118mm}{there is a continuous bounded function $c:[0,+\infty)  \to [0,+\infty)$ such that $$\|F(s,t,x)\| \le c(t)(1 + \|x\|_\alpha)\quad \mbox{ for } \  s\in[0,1], \ t\in [0,+\infty), \ x\in X^\alpha.$$}
\begin{definition}
We say that a continuous map $u:[0,+\infty) \to X^\alpha$ is \emph{a (global) mild solution} of the equation \ref{row-a-fs} starting at $x_0\in X^\alpha$, provided
\begin{equation*}
u(t) = S_A(t)x_0 + \int_0^t S_A(t - \tau)F(s,\tau,u(\tau)) \,d \tau \qquad\mathrm{for}\quad t\ge 0.
\end{equation*}
\end{definition}

In the following proposition we collect important facts concerning existence, continuity and compactness for the mild solutions of the equation \ref{row-a-fs}.
\begin{proposition}\label{th-exist1}
Under the above assumptions the following assertions hold. \\[-10pt]
\begin{itemize}
\item[(a)] For every $s\in[0,1]$ and $x\in X^\alpha$, the equation \ref{row-a-fs} admits a unique mild solution $u(\,\cdot\,;s, x):[0,+\infty)\to X^\alpha$ starting at $x$. \\[-7pt]
\item[(b)] If sequences $(x_n)$ in $X^\alpha$ and $(s_n)$ in $[0,1]$ are such that $x_n\to x_0$ in $X^\alpha$ and $s_n\to s_0$ when $n\to +\infty$, then given $t\ge 0$,
\begin{equation*}
u(t;s_n, x_n) \to u(t;s_0, x_0) \qquad\mathrm{as}\quad n\to +\infty,
\end{equation*}
and this convergence is uniform on the bounded subsets of $[0,+\infty)$. \\[-7pt]
\item[(c)] If $t > 0$ and $\Omega\subset X^\alpha$ is a bounded set, then $$\{u(t;s,x) \ | \ s\in[0,1], \ x\in\Omega\}$$ is a relatively compact subset of $X^\alpha$.
\end{itemize}
\end{proposition}
\begin{proof}
The proof of point $(a)$ is a consequence of \cite[Theorem 3.3.3]{Henry} and
\cite[Corollary 3.3.5]{Henry}. Points $(b)$ and $(c)$ are consequences of \cite[Proposition 4.1]{Kok2}.
\end{proof}

In the paper we will also use the following technical lemma.

\begin{lemma}\label{lem-conv}
Let $(w_n)$ in $C([0,+\infty),X)$ and $(x_n)$ in $X^\alpha$ be bounded sequences and let $u_n:[0,+\infty) \to X^\alpha$ be a map given by
$$u_n(t) := S_A(t) x_n + \int_0^t S_A(t - \tau)w_n(\tau)\,d \tau \qquad\mathrm{for}\quad t\in[0,+\infty).$$
Then the following assertions hold.
\begin{itemize}
\item[(i)]  The family of maps $\{u_n\}_{n\ge 1}$ is equicontinuous on $(0,+\infty)$.
\item[(ii)] If $(t_n)$ in $[a,+\infty)$, where $0 < a < +\infty$, is a sequence such that $t_n \to a$ as $n\to+\infty$, then the set $\{u_n(t_n) \ | \ n\ge 1\}$ is relatively compact in $X^\alpha$.
\end{itemize}
\end{lemma}
Before we start the proof of lemma, in the following proposition we collect useful properties of fractional spaces that will be applied in this paper.
\begin{proposition}\label{tw-prop-sect-oper}
The following assertions hold.
\begin{itemize}
\item[(a)] $S_A(t)X\subset X^{\alpha}$ for every $t>0$.\\[-9pt]
\item[(b)] If $x\in D(A^{\alpha})$ then $S_A(t)A^{\alpha}x = A^{\alpha}S_A(t)x$ for $t\ge 0$.\\[-9pt]
\item[(c)] There are $c>0$ and $M_\alpha > 0$ such that
$$A^{\alpha}S_A(t)\in \mathcal{L}(X) \ \mbox{ and } \ \|A^{\alpha}S_A(t)\|\le M_\alpha t^{-\alpha} e^{-ct} \qquad\mathrm{for}\quad t > 0.$$
\item[(d)] If $\beta_1, \beta_2\in \mathbb{R}$ and $\gamma := \max(\beta_1, \beta_2, \beta_1 + \beta_2)$, then $$A^{\beta_1 + \beta_2}x = A^{\beta_1}A^{\beta_2}x \qquad\mathrm{for}\quad x\in D(A^\gamma).$$
\end{itemize}
\end{proposition}
\begin{proof}
This is exactly \cite[Theorem 2.6.13]{Pazy} and \cite[Theorem 2.6.8]{Pazy}.
\end{proof}

\begin{proof}[Proof of Lemma \ref{lem-conv}.]
For the proof of point $(a)$ see \cite[Lemma 3.6]{Kok2}. To prove $(b)$ we show that $V:=\{A^\alpha u_n(t_n) \ | \ n\ge 1\}$ is a relatively compact subset of $X$. Let us take arbitrary $\varepsilon > 0$. By Proposition \ref{tw-prop-sect-oper} $(c)$ we have
\begin{equation}\label{afgh}
\left\|\int_{t'}^{t} A^\alpha S_A(t - \tau)w_n(\tau) \,d \tau \right\| \le \int_{t'}^{t} \frac{M_\alpha K }{(t - \tau)^\alpha} \,d \tau = (1 - \alpha)M_\alpha K (t - t')^{1 - \alpha}
\end{equation}
for $2a > t > t' \ge 0$, where $K = \sup\{\|w_n(\tau)\| \ | \ n\ge 1, \ \tau\in[0,2a]\}$. Hence there are $t_0\in(0,a)$ and $n_0 \ge 1$ such that
\begin{equation*}
\left\|\int_{t_0}^{t_n} A^\alpha S_A(t_n - \tau)w_n(\tau) \,d \tau \right\| \le \varepsilon \qquad\mathrm{for}\quad n\ge n_0.
\end{equation*}
Let us observe that, for any $n\ge 1$, one has
\begin{equation*}
\begin{split}
A^\alpha u_n(t_n) = & S_A(t_n)A^\alpha x_n + S_A(a - t_0)\left(\int_0^{t_0} A^\alpha S_A(t_n - a + t_0 - \tau)w_n(\tau) \,d \tau \right)\\
& + \int_{t_0}^{t_n} A^\alpha S_A(t_n - \tau)w_n(\tau)\,d \tau,
\end{split}
\end{equation*}
which, in view of \ref{afgh}, implies that
\begin{equation*}
\begin{aligned}
\{A^\alpha u_n(t_n) \ | \ n\ge n_0\}\subset S_A(t_0)D_1 + S_A(a - t_0)D_2 + B(0,\varepsilon),
\end{aligned}
\end{equation*}
where $B(0,\varepsilon):=\{x\in X \ | \ \|x\| \le \varepsilon\}$, $D_1:=\{S(t_n - t_0)A^\alpha x_n \ | \ n\ge 1\}$ and $$D_2:=\left\{S_A(t_n - a)\int_0^{t_0} A^\alpha S_A(t_0 - \tau)w_n(\tau) \,d \tau \ \Big| \ n\ge 1\right\}.$$
Since $(x_n)$ is bounded in $X^\alpha$, we see that $D_1$ is bounded in $X$. Furthermore \ref{afgh} implies that the set $D_2$ is bounded in $X$ as well. Hence the compactness of the semigroup and the fact that $\varepsilon > 0$ is arbitrary small imply that the set $V$ is relatively compact in $X$. As we noted at the beginning of the proof, this proves that $\{u_n(t_n)\}_{n\ge 1}$ is a relatively compact subset of $X^\alpha$ as desired.
\end{proof}

\section{Krasnosel'skii type degree formula}

In this section we consider the following differential equation
\begin{equation}\label{row-pert}
\dot u(t) = - A u(t) + F (u(t)), \qquad t > 0
\end{equation}
where $A:X\supset D(A)\to X$ is a positive sectorial operator with compact resolvents on a Banach space $X$ and $F:X^\alpha \to X$, where $\alpha\in (0,1)$, is a continuous map such that assumptions $(F1)$ and $(F2)$ are satisfied. From Proposition \ref{th-exist1} $(a)$, it follows that for every $x\in X^\alpha$, the equation \ref{row-pert} admits a unique mild solution $u(\,\cdot\,; x):[0,+\infty)\to X^\alpha$, such that $u(0; x) = x$. Let $\Phi_t:X^\alpha \to X^\alpha$ be {\em the Poincar\'e operator}, given for any $t > 0$, by $$\Phi_t(x) := u(t;x) \qquad\mathrm{for}\quad x\in X^\alpha.$$
Then Proposition \ref{th-exist1} $(b)$ and $(c)$ imply that $\Phi_t:X^\alpha \to X^\alpha$ is a completely continuous map for $t > 0$. After this comment, we are ready to provide the following Krasnosel'skii type degree formula for perturbations of sectorial operators.

\begin{theorem}\label{Th-kras-frac}
Let us assume that $U\subset X^\alpha$ is an open bounded set such that $0\notin (-A + F)(\partial U \cap D(A))$. Then there is  $\overline{t} > 0$ with the property that, if $t\in (0,\overline{t}]$, then $\Phi_t(x)\neq x$ for $x\in \partial U$ and
$$\mathrm{deg_{LS}}(I - \Phi_t, U) = \mathrm{deg_\alpha}(-A + F, U).$$
Here $\mathrm{deg_\alpha}$ is the topological degree for perturbations of sectorial operators. \footnote{\, For the definition and properties see Appendix.}
\end{theorem}

In the proof we will use the following version of Theorem \ref{Th-kras-frac}, where the nonlinear perturbation is defined only on the space $X$.

\begin{theorem}\label{th-cwiszewski2}
Let us consider the following differential equation
\begin{equation*}
\dot u(t) = - A u(t) + G (u(t)), \qquad t > 0
\end{equation*}
where $G:X \to X$ is a continuous map satisfying $(F1)$ and $(F2)$ (with $\alpha = 0$) and let us assume that $\Theta_t : X \to X$ is the associated Poincar\'e operator. If $U\subset X$ is a bounded open set such that $-Ax + G(x)\neq 0$ for $x\in \partial U \cap D(A)$, then there is $\overline{t} > 0$ with the property that, if $t\in (0,\overline{t}]$ then $\Theta_t(x)\neq x$ for $x\in \partial U$ and
\begin{equation*}
\mathrm{deg_{LS}}(I - \Theta_t, U) = \mathrm{deg_{LS}}(I - A^{-1}G, U).
\end{equation*}
\end{theorem}
\begin{proof}
The theorem is a combination of \cite[Theoem 4.5]{Cwiszewski-1} and \cite[Remark 4.8]{Cwiszewski-1}.
\end{proof}

In the proof of Theorem \ref{Th-kras-frac} we will also use the following {\em Volterra type inequality}.
\begin{lemma}\label{lem-volt-type-ineq}
Let us assume that $\alpha\in [0,1)$, $a\ge 0$, $b > 0$ and let $\phi\colon [0,T)\to [0,+\infty)$ be a continuous function such that
\begin{equation*}
\phi(t)\le a + b\int_{0}^t \frac{1}{(t - \tau)^\alpha} \phi(\tau) \,d \tau \qquad\mathrm{for}\quad  t\in (0, T).
\end{equation*}
Then there is a constant $K(\alpha,b,T)$ dependent from $\alpha$, $b$ and $T$ such that $$\phi(t)\le a K(\alpha,b,T) \qquad\mathrm{for}\quad t\in [0, T).$$
\end{lemma}
\begin{proof}
For the proof see \cite[Lemma 1.2.9]{DloCho}.
\end{proof}

\begin{proof}[Proof of Theorem \ref{Th-kras-frac}.] \emph{Step 1.} Let us observe that there is $\varepsilon_0 > 0$ such that
\begin{equation}\label{ja11}
-A x + S(\lambda\varepsilon)F(x)\neq 0 \qquad\mathrm{for}\quad \varepsilon \in (0,\varepsilon_0], \ \lambda\in [0,1], \ x\in \partial U\cap D(A).
\end{equation}
Suppose, contrary to our claim, that there are sequences $(\varepsilon_n)$ in $(0,1)$, $(\lambda_n)$ in $[0,1]$ and  $(x_n)$ in $\partial{U}\cap D(A)$ such that $\varepsilon_n \to 0^+$ as $n\to+\infty$ and
\begin{equation}\label{r-n-1}
-A x_n + S(\lambda_n\varepsilon_n)F(x_n) = 0 \qquad\mathrm{for}\quad n\ge 1.
\end{equation}
Since the operator $A$ is positive, by Proposition \ref{tw-prop-sect-oper} $(c)$, there is a constant $M > 0$ such that $\|S_A(t)\|\le M$ for $t\ge 0$. Then assumption $(F2)$ implies that $$\|x_n\|_1 = \|A x_n\| = \|S(\lambda_n\varepsilon_n)F(x_n)\| \le c M(1 + \|x_n\|_\alpha).$$
By the boundedness of $U\subset X^\alpha$, we infer that the sequence $(x_n)$ is bounded in $X^1$. Since the inclusion $X^1 \subset X^\alpha$ is compact (see \cite[Theorem 1.4.8]{Henry}), we deduce that there is $x_0\in \partial{U}$ such that $x_n \to x_0$ in $X^\alpha$. Writing the equation \ref{r-n-1} in the form
$$x_n - A^{-1}S(\lambda_n\varepsilon_n)F(x_n) = 0 \qquad\mathrm{for}\quad n\ge1,$$
and letting $n\to\infty$ one has $x_0 - A^{-1}F(x_0) = 0$, which contradicts the assumption. Consequently, by the homotopy invariance (see Theorem \ref{th-degal} $(D3)$),
\begin{equation}\label{dd3}
\mathrm{deg_\alpha}(-A + S(\varepsilon)F, U) = \mathrm{deg_\alpha}(-A + F, U) \qquad\mathrm{for}\quad \varepsilon\in(0,\varepsilon_0].
\end{equation}
\emph{Step 2.} Let us consider the following family of initial value problems
\begin{equation}\label{row-pert-2}
\left\{ \begin{aligned}
& \dot u(t) = - A u(t) + S(\lambda)F(u(t)), \qquad t > 0 \\
& u(0)=x \in X^\alpha. \end{aligned}\right.
\end{equation}
where $\lambda\in [0,1]$ is a parameter. Let us observe that the map $(\lambda,x)\mapsto S(\lambda)F(x)$, defined on $[0,1]\times X^\alpha$, satisfies assumptions $(F1)$ and $(F2)$. Therefore, by Proposition \ref{th-exist1} $(b)$ and $(c)$, the Poincar\'e operator  $\Psi_t:[0,1]\times \overline{U} \to X^\alpha$ for the equation
\ref{row-pert-2} is well-defined  and completely continuous for every $t\in (0,+\infty)$. Given $\varepsilon \in [0,1]$, let $\Psi_t^\varepsilon:[0,1]\times \overline{U} \to X^\alpha$ be defined by $\Psi_t^\varepsilon(\lambda,x):= \Psi_t(\lambda\varepsilon, x)$. We claim that there is $\varepsilon_1 \in (0, \varepsilon_0]$ and $t_0 > 0$ such that \begin{equation}\label{njnj}
\Psi_t^\varepsilon(\lambda, x)\neq x \qquad\mathrm{for}\quad \varepsilon\in (0, \varepsilon_1], \ t\in (0,t_0], \ \lambda\in [0,1] \text{ and } x\in \partial U.
\end{equation}
Indeed, suppose contrary to our claim that there are sequences $(\varepsilon_n)$, $(t_n)$ in $(0,1)$, $(\lambda_n)$ in $[0,1]$ and $(x_n)$ in $\partial U$ such that
$\varepsilon_n, t_n \to 0^+$, $\lambda_n \to \lambda_0$ as $n\to +\infty$ and
\begin{equation*}
\Psi_{t_n}^{\varepsilon_n}(\lambda_n, x_n) =  x_n \qquad\mathrm{for}\quad n\ge1.
\end{equation*}
Fix a real number $t_1\in (0,+\infty)$ and let $(k_n)$ be the sequence of integers given by $k_n := [t_1/t_n] + 1$ for $n\ge 1$. It is not difficult to see that $k_n t_n \ge t_1$ and $k_n t_n \to t_1$ as $n \to \infty$. If we put $r_n := k_n t_n - t_1$, then
\begin{equation*}
x_n =  \Psi_{k_n t_n}(\varepsilon_n\lambda_n, x_n) = \Psi_{t_1}(\varepsilon_n\lambda_n, \Psi_{r_n}(\varepsilon_n\lambda_n, x_n))\qquad\mathrm{for}\quad n\ge 1.
\end{equation*}
In view of assumption $(F2)$ and Lemma \ref{lem-volt-type-ineq}, the set $\{\Psi_{r_n}(\varepsilon_n\lambda_n, x_n)\}_{n\ge 1}\subset X^\alpha$ is bounded, which along with the compactness of $\Psi_{t_1}$ implies that the sequence $(x_n)$ is relatively compact in $X^\alpha$. Hence, there is a subsequence $(x_{n_l})$ of $(x_n)$ such that $x_{n_l} \to x_0\in \partial U$ as $l \to \infty$. Let us take an arbitrary $t\in (0,+\infty)$ and let $(m_l)$ be the sequence of integers given by $m_l:=[t/t_{n_l}] + 1$. Then $m_l t_{n_l} \ge t$ and $m_l t_{n_l} \to t$ as $l \to \infty$. Furthermore
\begin{equation}\label{r-n-2}
\Psi_{m_l t_{n_l}}(\varepsilon_{n_l}\lambda_{n_l}, x_{n_l}) = x_{n_l} \qquad\mathrm{for}\quad l\ge 1.
\end{equation}
Using Proposition \ref{th-exist1} $(b)$ and letting $l\to \infty$ in \ref{r-n-2} one has
$x_0 = \Psi_{t}(0,x_0)$. Since $t\in (0,+\infty)$ is arbitrary, it implies that
\begin{equation*}
x_0 = S_A(t)x_0 + \int_0^t S_A(t - \tau)F(x_0)\,d \tau \qquad\mathrm{for}\quad t\in (0,+\infty),
\end{equation*}
and, in particular,
\begin{equation*}
\frac{x_0 - S_A(t)x_0}{t} = \frac{1}{t}\int_0^t S_A(t - \tau)F(x_0)\,d \tau \qquad\mathrm{for}\quad t\in (0,+\infty).
\end{equation*}
Letting $t \to 0^+$ yields
$x_0 \in \partial U\cap D(A)$ and $-Ax_0 + F(x_0) = 0$, which is a contradiction and the claim follows. Therefore using \ref{njnj} together with homotopy invariance of topological degree, one has
\begin{equation}\label{row-deg-form-1}
\mathrm{deg_{LS}}(I - \Phi_t, U) = \mathrm{deg_{LS}}(I - \Psi_t^\varepsilon(0,\,\cdot\,), U) = \mathrm{deg_{LS}}(I - \Psi_t^\varepsilon(1,\,\cdot\,), U),
\end{equation}
for every $\varepsilon\in (0,\varepsilon_1]$ and $t \in (0, t_0]$. \\[5pt]
\noindent\emph{Step 3.} Let us take $\varepsilon \in (0,\varepsilon_1]$ and define
$u_\varepsilon (t;x) := \Psi_t^\varepsilon(1,x)$ for $t\in [0,+\infty)$ and $x\in X^\alpha$. Then
$$u_\varepsilon(t;x) = S_A(t)x + \int_0^t S_A(t - \tau)S(\varepsilon)F(u_\varepsilon(\tau;x))\,d \tau
\qquad\mathrm{for}\quad t\in [0,+\infty),$$ which together with Proposition \ref{tw-prop-sect-oper} $(a)$ and $(b)$ gives
\begin{equation}\label{row-var-const-1}
A^\alpha u_\varepsilon(t;x) = S_A(t)A^\alpha x + \int_0^t S_A(t - \tau)A^\alpha S(\varepsilon)F(A^{-\alpha}A^\alpha u_\varepsilon(\tau;x))\,d \tau
\end{equation}
for $t\in [0,+\infty)$. Let us consider the following initial value problem
\begin{equation}\label{row-pert-3}
\left\{ \begin{aligned}
& \dot u(t) = - A u(t) + G(u(t)), \qquad t > 0\\
& u(0)=x\in X, \end{aligned} \right.
\end{equation}
where $G:X\to X$ is given by
$$G(x) := A^\alpha S(\varepsilon)F(A^{-\alpha}x) \qquad\mathrm{for}\quad x\in X.$$
It is not difficult to see that $G$ is continuous and satisfies assumptions $(G1)$ and $(G2)$ (with $\alpha = 0$). Indeed, if $x_0\in X$ then there is a neighborhood $V\subset X^\alpha$ of $A^{-\alpha}x_0$ and $L > 0$ such that
\begin{equation*}
\|F(x) - F(y)\|\le L\|x - y\|_\alpha \qquad\mathrm{for}\quad x,y\in V.
\end{equation*}
Proposition \ref{tw-prop-sect-oper} $(c)$ implies that $A^\alpha S(\varepsilon)\in \mathcal{L}(X)$ and there is $M_\alpha > 0$ such that $$\|A^\alpha S(\varepsilon)\| \le M_\alpha \varepsilon^{-\alpha}.$$ It is not difficult to see that $A^\alpha V$ is a neighborhood of $x_0$ and, for any $x,y\in A^\alpha V$,
\begin{align*}
\|G(x) - G(y)\| & = \|A^\alpha S(\varepsilon)F(A^{-\alpha}x) - A^\alpha S(\varepsilon)F(A^{-\alpha}y)\| \\
& \le M_\alpha \varepsilon^{-\alpha}\|F(A^{-\alpha}x) - F(A^{-\alpha}y)\|  \\
& \le L M_\alpha \varepsilon^{-\alpha}\|A^{-\alpha}x - A^{-\alpha}y\|_\alpha \\
& = L M_\alpha \varepsilon^{-\alpha}\|x - y\|,
\end{align*}
which proves that $G$ is locally Lipschitz. Furthermore, using assumption $(F2)$, gives
\begin{align*}
\|G(x)\| & = \|A^\alpha S(\varepsilon)F(A^{-\alpha}x)\| \le M_\alpha \varepsilon^{-\alpha}\|F(A^{-\alpha}x)\| \\
& \le c M_\alpha \varepsilon^{-\alpha}(1 + \|A^{-\alpha}x\|_\alpha) = c M_\alpha \varepsilon^{-\alpha}(1 + \|x\|),
\end{align*}
for any $x\in X$ and consequently $G$ satisfies $(F2)$. Let $\Theta_t:X \to X$ be the Poincar\'e operator associated with the equation \ref{row-pert-3}.

In view of the equation \ref{row-var-const-1} and the uniqueness of solutions for \ref{row-pert-3} (see Proposition \ref{th-exist1} $(a)$), we infer that
\begin{equation}\label{equ-phi-theta}
A^\alpha\Psi_t^\varepsilon (1,A^{-\alpha}x) = \Theta_t(x) \qquad\mathrm{for}\quad x\in X,\quad t\in [0,+\infty).
\end{equation}
Let us note that $-A x + G(x) \neq 0$ for $x\in \partial(A^\alpha U)\cap D(A)$. Otherwise there would be $x_0\in \partial(A^\alpha U)\cap D(A)$ such that $-A x_0+ G(x_0) = 0$ and hence
\begin{equation}\label{fra}
-A x_0 + A^\alpha S(\varepsilon)F(A^{-\alpha} x_0) = 0.
\end{equation}
On the other hand, Theorem \ref{tw-prop-sect-oper} $(d)$ gives
$A^\alpha A^{-1} x = A^{\alpha - 1} x = A^{-1}A^\alpha x$ for $x\in D(A)$, which along with \ref{fra} implies that $y_0 = A^{-\alpha} x_0 \in\partial U \cap D(A)$ and
\begin{equation*}
-A y_0 + S(\varepsilon)F(y_0) = 0,
\end{equation*}
contrary to \ref{ja11}. Therefore, assumptions of Theorem \ref{th-cwiszewski2} are satisfied and there is  $\overline{t} \in (0, t_0]$ with the property that, if $t\in (0,\overline{t}]$ then $\Theta_t(x)\neq x$ for $x\in \partial(A^\alpha U)$ and
\begin{equation}\label{deg1a}
\mathrm{deg_{LS}}(I - \Theta_t, A^\alpha U) = \mathrm{deg_{LS}}(I - A^{-1}G, A^\alpha U).
\end{equation}
In view of \ref{equ-phi-theta} this implies that $\Psi_t^\varepsilon (1,x)\neq x$ for $t\in (0,\overline{t}]$, $x\in \partial U$ and
\begin{equation}\label{deg2}
\mathrm{deg_{LS}}(I - \Psi_t^\varepsilon(1,\,\cdot\,), U) = \mathrm{deg_{LS}}(I - \Theta_t, A^\alpha U).
\end{equation}
Combining \ref{row-deg-form-1}, \ref{deg2} and \ref{deg1a}, we infer that, for any $t\in(0,\overline{t}]$
\begin{equation}\label{deg3}
\begin{aligned}
\mathrm{deg_{LS}}(I - \Phi_t, U) & = \mathrm{deg_{LS}}(I - A^{-1}G, A^\alpha U) \\ & = \mathrm{deg_{LS}}(I - A^{\alpha - 1} S(\varepsilon)F(A^{-\alpha} \ \cdot \ ),A^\alpha U).
\end{aligned}
\end{equation}
On the other hand, from \ref{dd3} it follows that
\begin{align*}
\mathrm{deg_{LS}}(I - A^{\alpha - 1} S(\varepsilon)F(A^{-\alpha} \ \cdot \ ),A^\alpha U) & = \mathrm{deg_{LS}}(I - A^{-1}S(\varepsilon)F, U) \\
& = \mathrm{deg_\alpha}(-A + S(\varepsilon)F, U) \\
& = \mathrm{deg_\alpha}(-A + F, U),
\end{align*}
for $\varepsilon\in(0,\varepsilon_1]$, which together with \ref{deg3} imply that $$\mathrm{deg_{LS}}(I - \Phi_t, U) = \mathrm{deg_\alpha}(-A + F, U) \qquad\mathrm{for}\quad t\in (0,\overline{t}]$$ and the proof is completed.
\end{proof}

\section{Averaging principle for perturbations of sectorial operators}

Let us consider the following family of differential equations
\begin{equation}\label{glowne}
\dot u(t) = - \lambda A u(t) + \lambda F (t,u(t)), \qquad t > 0
\end{equation}
where $\lambda\in[0,1]$ is a parameter, $A:X\supset D(A) \to X$ is a positive sectorial operator with compact resolvents and $F:[0,+\infty)\times X^\alpha \to X$, $\alpha\in(0,1)$, is a continuous map satisfying $(F1)$, $(F2)$ and the following assumption \\[5pt]
\noindent\makebox[9mm][l]{$(F3)$}\parbox[t]{118mm}{there is $T > 0$ such that $F(t + T,x) = F(t, x)$  for $t\in[0,+\infty)$,  $x\in X^\alpha$.}\\[5pt]
It can be easily seen that for any $\lambda\in (0,1]$, the operator $\lambda A$ is also sectorial and $$\Re\sigma(\lambda A) = \lambda\Re\sigma(A) > 0.$$
By Proposition \ref{th-exist1} $(a)$, for any $\lambda\in(0,1]$ and $x\in X^\alpha$ there is a unique mild solution $u(\,\cdot\, ; \lambda, x):[0,+\infty) \to X^\alpha$ of the equation \ref{glowne} such that $u(0 ; \lambda, x) = x$. Given $\lambda\in (0,1]$ and $t\in[0,+\infty)$, define the Poincar\'e operator
$\Phi^\lambda_t:X^\alpha \to X^\alpha$ by $$\Phi^\lambda_t(x) := u(t ;\lambda,x) \qquad\mathrm{for}\quad x\in X^\alpha.$$
From Proposition \ref{th-exist1} $(b)$ and $(c)$ it follows easily that $\Phi^\lambda_T$ is completely continuous.
Now we are ready to prove the following {\em the averaging principle}, which is an effective topological tool to study the existence of $T$-periodic solutions.

\begin{theorem}\label{th-cont}{\em (Averaging principle)}
If $U\subset X^\alpha$ is an open bounded set such that $0 \notin (-A + \widehat{F})(\partial U \cap D(A))$, then  there is $\lambda_0 > 0$ such that if $\lambda\in (0,\lambda_0]$ then $\Phi^\lambda_T(x) \neq x$ for $x\in\partial U$ and $$\mathrm{deg_{LS}}(I - \Phi^\lambda_T, U) = \mathrm{deg_\alpha}(-A + \widehat{F}, U).$$
\end{theorem}

\begin{proof} \emph{Step 1.} Let us consider the following family of equations
\begin{equation}\label{row-pert-4}
\dot v(t) = - \lambda A v(t) + \lambda\widetilde{F}(\mu,t,v(t)) \qquad\mbox { na } \quad [0,+\infty)
\end{equation}
where $\lambda\in[0,1]$ is a parameter and $\widetilde{F}\colon [0,1]\times [0,+\infty)\times X^\alpha \to X$ is given by
$$\widetilde{F}(\mu,t,x) = \mu F (t,x) + (1 - \mu)\widehat{F}(x) \qquad\mathrm{for}\quad \mu\in[0,1], \ t\in[0,+\infty), \ x\in X^\alpha.$$
It is not difficult to see that $\widetilde{F}$ satisfies assumptions $(F1)$ and $(F2)$. Indeed, if $x_0\in X^\alpha$ then there are a neighborhood $V\subset X^\alpha$ of $x_0$ and $L > 0$ such that
\begin{equation*}
\|F(t,x) - F(t,y)\| \le L \|x - y\|_\alpha \qquad\mathrm{for}\quad t\in [0,T],\quad x,y\in V.
\end{equation*}
It follows that, for any $\mu\in [0,1]$ and $x, y\in V$, one has
\begin{align*}
\|\widetilde{F}(\mu,t,x) - \widetilde{F}(\mu,t,y)\| & =
\|\mu F (t,x) - \mu F (t,y) + (1 - \mu)\widehat{F}(x) - (1 - \mu)\widehat{F}(y)\| \\
& \le \mu\|F(t,x) - F(t,y)\| + (1 - \mu)\|\widehat{F}(x) - \widehat{F}(y)\| \\
& \le L \|x - y\|_\alpha + \frac{1}{T}\int_0^T \|F(\tau,x) - F(\tau,y)\| \,d \tau \\
& \le 2 L \|x - y\|_\alpha,
\end{align*}
which proves that $\widetilde{F}$ satisfies assumption $(F1)$. In order to show that $(F2)$ holds, observe that for any $t\in [0,+\infty)$, $\mu\in [0,1]$ and $x\in X^\alpha$
\begin{align*}
\|\widetilde{F}(\mu,t,x)\| & =
\|\mu F (t,x) + (1 - \mu)\widehat{F}(x)\| \le \mu\|F(t,x)\| + (1 - \mu)\|\widehat{F}(x)\| \\
& \le c(t)(1 + \|x\|_\alpha) + \frac{1}{T}\int_0^T c(\tau)(1 + \|x\|_\alpha) \,d \tau \\
& = c_0(t)(1 + \|x\|_\alpha),
\end{align*}
where $c_0(t) := \left(c(t) + \frac{1}{T}\int_0^T c(\tau) \,d \tau\right)$ for $t\in [0,+\infty)$. This proves that $(F2)$ holds. \\[5pt]
\emph{Step 2.} Let us observe that from Step 1 and Proposition \ref{th-exist1} $(a)$ it follows that, for any $x\in X^\alpha$, $\lambda\in[0,1]$ and $\mu\in[0,1]$, there is a unique mild solution $v(\,\cdot \, ;\lambda,\mu,x):[0,+\infty) \to X^\alpha$ of the equation \ref{row-pert-4} such that $v(0 ;\lambda,\mu,x) = x$. If $\Psi^\lambda_T:[0,1]\times X^\alpha \to X^\alpha$ is the associated Poincar\'e operator, then Proposition \ref{th-exist1} $(b)$ and $(c)$ prove that $\Psi^\lambda_T$ is a completely continuous map.
We show that there is $\lambda_1 > 0$ such that
\begin{equation*}
\Psi^\lambda_T(\mu, x) \neq x \qquad\mathrm{for}\quad \lambda\in (0,\lambda_1], \ \mu\in [0,1] \ \text{ and } \ x\in\partial U.
\end{equation*}
Suppose, contrary to our claim, that there are sequences $(\lambda_n)$ in $(0,1)$, $(\mu_n)$ in $[0,1]$, $(x_n)$ in $\partial U$ such that $\lambda_n \to 0^+$, $\mu_n \to \mu_0$ as $n\to +\infty$ and
\begin{equation}\label{r-l-4}
\Psi^{\lambda_n}_T(\mu_n, x_n) = x_n \qquad\mathrm{for}\quad n\ge 1.
\end{equation}
Given $n\ge 1$, let us define $v_n: = v(\,\cdot\, ;\lambda_n,\mu_n,x_n)$ and observe that the set
\begin{equation}\label{zbior}
E:=\{\widetilde{F}(\mu_n,\tau,v_n(\tau)) \ | \ \tau\in [0,T], \ n\ge 1\}
\end{equation}
is bounded in $X$. Indeed, by Proposition \ref{tw-prop-sect-oper} $(c)$, there are $M,M_\alpha > 0$ such that
\begin{equation}\label{esti1}
\|S_A(t)\| \le M, \qquad \|A^\alpha S_A(t)\| \le M_\alpha/t^\alpha  \qquad\mathrm{for}\quad t > 0.
\end{equation}
If $R > 0$ is such that $\|x_n\|_\alpha\le R$ for $n\ge 1$, then
\begin{align*}
\|v_n(t)\|_\alpha & \le  \|S_A(\lambda_n t)A^\alpha x_n\|
+ \lambda_n\int_0^t \|A^\alpha S_A(\lambda_n(t - \tau))\|\|\widetilde{F}(\mu_n,\tau,v_n(\tau))\|\,d \tau \\
& \le  M \|x_n\|_\alpha
+ \lambda_n^{1 - \alpha}\int_0^t \frac{M_\alpha}{(t - \tau)^\alpha}\|\widetilde{F}(\mu_n,\tau,v_n(\tau))\|\,d \tau \\
& \le M R
+ \int_0^{t} \frac{M_\alpha}{(t - \tau)^\alpha}c_0(\tau)(1 + \|v_n(\tau)\|_\alpha) \,d \tau \\
& \le M R + \frac{K M_\alpha}{1-\alpha} T^{1-\alpha}
+ \int_0^{t} \frac{K M_\alpha }{(t - \tau)^\alpha}\|v_n(\tau)\|_\alpha \,d \tau,
\end{align*}
where $K := \sup_{\tau\in [0,T]}c_0(\tau)$. Hence, by Lemma \ref{lem-volt-type-ineq}, we infer that there is $C > 0$  such that $\|v_n(t)\|_\alpha \le C$ for $t\in [0,T]$ and $n\ge 1$. Consequently
\begin{align*}
\|\widetilde{F}(\mu_n,\tau,v_n(\tau))\| \le c_0(\tau)(1 + \|v_n(\tau)\|_\alpha) \le K(1 + C) \ \text{ for } \
s\in [0,T], \ n\ge 1,
\end{align*}
which implies that the set \ref{zbior} is bounded. Let us write \ref{r-l-4} in the following form
\begin{equation}\label{r-l-6}
x_n = v_n(T) = S_A(\lambda_n T)x_n +
\lambda_n\int_0^T S_A(\lambda_n(T - \tau))\widetilde{F}(\mu_n,\tau,v_n(\tau)) \,d \tau.
\end{equation}
We claim that
\begin{equation}\label{rr4}
v_n(t) = v_n(t + T) \qquad\mathrm{for}\quad t\in[0,+\infty).
\end{equation}
Indeed, using the fact that $\widetilde{F}$ is $T$-periodic, one has
\begin{align*}
v_n(t + T) & = S_A(\lambda_n t)v_n(T) +
\lambda_n\int_T^{t + T} S_A(\lambda_n(t+T -\tau))\widetilde{F}(\mu_n,\tau,v_n(\tau)) \,d \tau \\
& = S_A(\lambda_n t)x_n +
\lambda_n\int_0^t S_A(\lambda_n(t - \tau))\widetilde{F}(\mu_n,\tau + T,v_n(\tau + T)) \,d \tau \\
& = S_A(\lambda_n t)x_n +
\lambda_n\int_0^t S_A(\lambda_n(t - \tau))\widetilde{F}(\mu_n,\tau,v_n(\tau + T)) \,d \tau.
\end{align*}
This implies that the map $v_n( \, \cdot \, + T):[0,+\infty) \to X^\alpha$ is a mild solution of the equation \ref{row-pert-4} starting at $x_n$ and consequently \ref{rr4} follows from Proposition \ref{th-exist1} $(a)$.
Therefore, combining \ref{r-l-4} and \ref{rr4} yields $x_n = v_n(T) = v_n(kT)$, which by the Duhamel formula gives
\begin{equation}\label{r-l-8}
x_n = S_A(\lambda_n k T)x_n +
\lambda_n\int_0^{kT} S_A(\lambda_n(kT - \tau))\widetilde{F}(\mu_n,\tau,v_n(\tau)) \,d \tau.
\end{equation}
Let $(k_n)$ be the sequence of integers given by
$$k_n := [1/(2 \lambda_n)] + 1 \qquad\mathrm{for}\quad n\ge 1.$$
It is not difficult to see that the sequence $t_n := k_n\lambda_n T$ is such that $t_n \ge T/2$ and $t_n \to T/2$ as $n\to +\infty$. Let us write the formula \ref{r-l-8} in the following form
\begin{equation}
\begin{aligned}\label{rr9}
x_n  & = S_A(\lambda_n k_n T) x_n +
\int_0^{\lambda_n k_n T} S_A(\lambda_n k_nT - \tau) \widetilde{F}(\mu_n,\tau/\lambda_n,v_n(\tau/\lambda_n)) \,d \tau \\
& = S_A(t_n)x_n +
\int_0^{t_n} S_A(t_n - \tau)\widetilde{F}(\mu_n,\tau/\lambda_n,v_n(\tau/\lambda_n)) \,d \tau.
\end{aligned}
\end{equation}
Let us consider the sequence of maps $w_n(\tau) := \widetilde{F}(\mu_n,\tau/\lambda_n,v_n(\tau/\lambda_n))$ for $\tau\ge 0$.
From \ref{rr4} and assumption $(F3)$ we infer that
\begin{align*}
\{w_n(\tau) \ | \ n\ge 1, \ \tau\in[0,+\infty)\} & = \{\widetilde{F}(\mu_n,\tau/\lambda_n,v_n(\tau/\lambda_n)) \ | \ \tau\in [0,+\infty), \ n\ge 1\} \\
& \subset \{\widetilde{F}(\mu_n,\tau,v_n(\tau)) \ | \ \tau\in [0,T], \ n\ge 1\} = E.
\end{align*}
As it was proved in \ref{zbior}, the set $E$ is bounded in $X$, which implies that the sequence $(w_n)$ is bounded in the space $C([0,+\infty),X)$. In view of \ref{rr9} and Lemma \ref{lem-conv} $(ii)$ applied for $a:=T/2$, we infer that the sequence $(x_n)$, contained in $\partial U$, is relatively compact in $X^\alpha$. Therefore, without loss of generality, we can suppose that there is $x_0\in \partial U$ such that $x_n\to x_0$ as $n\to+\infty$.
Let us take arbitrary $t\in [0,T]$. By \ref{esti1}
\begin{align*}
\|w_n(t) - x_n\|_\alpha & \le \|S_A(\lambda_n t)A^\alpha x_n - A^\alpha x_n\| + \lambda_n\int_0^t \|A^\alpha S_A(\lambda_n t - \lambda_n \tau)\widetilde{F}(\mu_n,\tau,w_n(\tau)) \| \,d \tau \\
& \le \|S_A(\lambda_n t)A^\alpha x_n - A^\alpha x_n\|
+ \lambda_n^{1 - \alpha}\int_0^t M_\alpha (t - \tau)^{-\alpha}\|\widetilde{F}(\mu_n,\tau,w_n(\tau))\| \,d \tau \\
& \le \|S_A(\lambda_n t)A^\alpha x_n - A^\alpha x_n\| + \lambda_n^{1 - \alpha}\int_0^t K M_\alpha (t - \tau)^{-\alpha} \,d \tau \\
& \le \|S_A(\lambda_n t) A^\alpha x_n - A^\alpha x_n\|
+ \frac{K M_\alpha}{1 - \alpha}\lambda_n^{1 - \alpha}T^{1 - \alpha},
\end{align*}
where $K := \sup\{\|\widetilde{F}(\mu_n,\tau,w_n(\tau))\| \ | \ \tau\in[0,T], \ n\ge 1\}$. This in turn, implies that $(w_n)$ converges in uniformly on $[0,T]$ to a constant function identically equal to $x_0$. If we write \ref{r-l-6} in the following form
\begin{equation*}
\frac{x_n - S_A(\lambda_n T)x_n}{\lambda_n T} =
\frac{1}{T}\int_0^T S_A(\lambda_n(T - \tau))\widetilde{F}(\mu_n,\tau,w_n(\tau)) \,d \tau,
\end{equation*}
then, by the standard properties of $C_0$ semigroups
\begin{align}\label{r-l-11}
A\left(\frac{\int_0^{\lambda_n T}S_A(\tau)x_n \,d \tau}{\lambda_n T} \right) =
\frac{1}{T}\int_0^T S_A(\lambda_n(T - \tau))\widetilde{F}(\mu_n,\tau,w_n(\tau)) \,d \tau.
\end{align}
It is not difficult to see that
\begin{equation}\label{r-l-12}
\frac{\int_0^{\lambda_n T}S_A(\tau)x_n \,d \tau}{\lambda_n T} \to x_0 \qquad\text{as}\quad n\to \infty.
\end{equation}
Since $w_n(t)\to x_0$ uniformly on $[0,T]$, the sequence of maps
\begin{equation*}
\tau \mapsto S_A(\lambda_n(T - \tau))\widetilde{F}(\mu_n,\tau,w_n(\tau))
\end{equation*}
converges uniformly on $[0,T]$ to the map $\tau \mapsto \widetilde{F}(\mu_0,\tau,x_0)$. Therefore
\begin{equation}\label{r-l-14}
\frac{1}{T}\int_0^T S_A(\lambda_n(T - \tau))\widetilde{F}(\mu_n,\tau,u_n(\tau)) \,d \tau \to \frac{1}{T}\int_0^T \widetilde{F}(\mu_0,\tau,x_0) \,d \tau = \widehat{F}(x_0)
\end{equation}
as $n\to \infty$. Combining \ref{r-l-11}, \ref{r-l-12}, \ref{r-l-14}  and using the fact that $A$ is a closed operator, we infer that $x_0\in D(A)$ and $-A x_0 + \widehat{F}(x_0) = 0$. This is a contradiction because $x_0\in \partial U\cap D(A)$ and the assertion of Step 2 is proved. \\[5pt]
\emph{Step 3.} By Step 2 and the homotopy invariance of topological degree
\begin{equation}\label{r-l-15}
\mathrm{deg_{LS}}(I - \Phi^\lambda_T, U) = \mathrm{deg_{LS}}(I - \Psi^\lambda_T(1,\,\cdot\,), U) = \mathrm{deg_{LS}}(I - \Psi^\lambda_T(0,\,\cdot\,), U)
\end{equation}
for $\lambda\in (0,\lambda_1]$. Since $\Psi^\lambda_T(0,\,\cdot\,)$ is the Poincar\'e operator for the autonomous equation \ref{row-pert-4} with $\mu=0$, one has
\begin{equation}\label{r-l-16}
\Psi^\lambda_T(0,x) = \Psi^1_{\lambda T}(0,x) \qquad\mathrm{for}\quad x\in X^\alpha.
\end{equation}
Hence, Theorem \ref{Th-kras-frac} implies the existence of $\lambda_0\in (0,\lambda_1]$ such that, if $\lambda\in (0,\lambda_0]$ then $\Psi^1_{\lambda T}(0,x) \neq x$ for $x\in\partial U$ and
\begin{equation}\label{r-l-17}
\mathrm{deg_{LS}}(\Psi^1_{\lambda T}(0,\cdot), U) = \mathrm{deg_\alpha}(-A + \widehat{F}, U).
\end{equation}
Combining \ref{r-l-15}, \ref{r-l-16} and \ref{r-l-17}, we deduce that, if $\lambda\in (0,\lambda_0]$ then
\begin{align*}
\mathrm{deg_{LS}}(I - \Phi^\lambda_T, U) & = \mathrm{deg_{LS}}(I - \Psi^\lambda_T(0,\cdot), U)
= \mathrm{deg_{LS}}(I - \deg(\Psi^1_{\lambda T}(0,\cdot), U)) \\
& = \mathrm{deg_\alpha}(-A + \widehat{F}, U)
\end{align*}
and the proof of theorem is completed.
\end{proof}

\section{Global continuation principle}

In this section we will consider differential equation of the following form
\begin{equation}\label{nowe}
    \dot u(t) = - A u(t) + F (t,u(t)), \quad t > 0.
\end{equation}
where $A:X\supset D(A) \to X$ is a positive sectorial operator with compact resolvent and $F:[0,+\infty)\times X^\alpha \to X$, $\alpha\in(0,1)$, is a continuous map satisfying assumptions $(F1)-(F3)$.
Using the averaging principle for the perturbations of sectorial operators, that was obtain in the previous section, we derive the following {\em global continuation principle}, that is an effective criterion for the existence of $T$-periodic solutions for the equation \ref{nowe}.
\begin{theorem}{\em (Global continuation principle)}\label{th-cont2}
Let us consider the family of bounded linear operators $\{F_\infty(t): X^\alpha \to X\}_{t\ge 0}$ such that the map $t\mapsto F_\infty(t) \in \mathcal{L}(X^\alpha,X)$ is $T$-periodic, continuous on $[0,+\infty)$ and
\begin{equation}\label{nnrr}
\lim_{\|x\|_\alpha\to +\infty} \frac{\|F(t,x) - F_\infty(t)x\|}{\|x\|_\alpha} = 0 \qquad \text{uniformly for } \ t\in[0,T].
\end{equation}
If the parameterized linear problem
\begin{equation}\label{nnrr2}
\left\{\begin{aligned}
\dot u(t) & = \lambda(A + F_\infty(t))u(t), \qquad t\in[0,+\infty) \\
u(0) & = u(T)
\end{aligned}\right.
\end{equation}
does not admit nontrivial solution for any $\lambda\in(0,1]$ and $\mathrm{Ker}(\widehat A + \widehat F_\infty) = \{0\}$, then the equation \ref{nowe} admits a $T$-periodic mild solution.
\end{theorem}
\begin{proof} \emph{Step 1. }Let us start with the claim that there is $R_0 > 0$ such, that for any $\lambda\in(0,1]$, the equation \ref{glowne} does not admit $T$-periodic solutions starting from the set $\{x\in X^\alpha \ | \ \|x\|_\alpha \ge R_0\}$. Otherwise there would be sequences $(\lambda_n)$ in $(0,1]$ and $(u_n)$ in $X^\alpha$ such that $u_n$ is a $T$-periodic solution of \ref{glowne} and $\|u_n\|_\alpha\to +\infty$ as $n\to +\infty$. Let us define the sequence of maps $F_n:[0,+\infty)\times X^\alpha \to X$ by
$$F_n(t,x):= \|u_n\|^{-1}_\infty F(t, \|u_n\|_\infty x) \qquad\mathrm{for}\quad t\in[0,+\infty), \ x\in X^\alpha,$$
where $\|u_n\|_\infty:= \sup_{t\in [0,+\infty)}\|u_n(t)\|_\alpha$. Then, by the Duhamel formula,
\begin{equation}\label{roronn2}
v_n(t) = S_A(\lambda_n t)x_n + \lambda_n\int_0^tS_A(\lambda_n(t - \tau))F_n(\tau,v_n(\tau)) \, d\tau
\end{equation}
where $v_n(t) := u_n(t)/\|u_n\|_\infty$ and $x_n:= u_n(0)/\|u_n\|_\infty$. By assumption $(F2)$,
\begin{equation}\label{roroee}
\begin{aligned}
\|F_n(t,v_n(t))\| & = \|\|u_n\|^{-1}_\infty F(t, \|u_n\|_\infty v_n(t))\| \\
& \le \|u_n\|^{-1}_\infty c(t)(1 + \|u_n\|_\infty \|v_n(t)\|_\alpha) \\
& \le K_1 K_2 + K_1 := K \qquad\mathrm{for}\quad n\ge 1, \ t\in[0,+\infty),
\end{aligned}
\end{equation}
where $K_1:=\sup_{\tau\in[0,+\infty)}c(\tau)$ and $K_2:=\sup_{n\ge 1} \|u_n\|^{-1}_\infty$. Passing if necessary to a subsequence, we can assume that there is $\lambda_0\in[0,1]$ such that $\lambda_n\to\lambda_0$ as $n\to+\infty$. \\[5pt]
\emph{Case A.} Suppose that $\lambda_0\in(0,1]$. From the equation \ref{roronn2}, it follows that
$$v_n(t) = S_A(\lambda_n t)x_n + \int_0^{t \lambda_n}S_A(\lambda_nt - \tau)F_n(\tau/\lambda_n,v_n(\tau/\lambda_n)) \, d\tau \qquad\mathrm{for}\quad t\ge 0$$
which together with Lemma \ref{lem-conv} $(i)$ and \ref{roroee} imply that the family $\{v_n(\,\cdot\,/\lambda_n)\}_{n\ge 1}$ is equicontinuous on $(0,+\infty)$. Since $(\lambda_n)$ is contained in $(0,1]$ and $\lambda_0\neq 0$, we infer that the family $\{v_n\}_{n\ge 1}$ is also equicontinuous at each point of the interval $(0,+\infty)$. Furthermore, given $t\in(0,+\infty)$, we use \ref{roroee} and we apply Lemma \ref{lem-conv} $(ii)$ with $t_n := \lambda_n t$ to obtain relative compactness of the sequence $(v_n(t))$ in $X^\alpha$. Since $v_n$ is $T$-periodic, the sequence of maps $(v_n)$, restricted to the interval $[0,T]$, satisfies assumptions of Ascoli-Arzela Theorem and therefore contains a subsequence that is convergent in $C([0,T],X^\alpha)$. Without loss of generality we can assume that there is continuous $T$-periodic $v_0:[0,+\infty) \to X^\alpha$ such that $v_n\to v_0$ uniformly on $[0,+\infty)$. Since $\|v_n\|_\infty = 1$ for $n\ge 1$, one has $\|v_0\|_\infty = 1$. Let us observe that, by the equation \ref{roronn2}
\begin{equation}\label{efefeff}
\begin{aligned}
v_n(t) & = S_A(\lambda_n t)x_n + \lambda_n\int_0^tS_A(\lambda_n(t - \tau))(F_n(\tau,v_n(\tau)) - F_\infty(\tau)v_n(\tau))\, d\tau \\
& \qquad + \lambda_n\int_0^tS_A(\lambda_n(t - \tau))F_\infty(\tau)v_n(\tau)\, d\tau \qquad\mathrm{for}\quad t\ge 0.
\end{aligned}
\end{equation}
From \ref{nnrr}, we infer that, given $\varepsilon > 0$ there is $m_\varepsilon > 0$ such that $\|F(t,x) - F_\infty(t)x\| \le \varepsilon\|x\|_\alpha + m_\varepsilon$ for $t\in[0,+\infty)$ and $x\in X^\alpha$. Therefore
\begin{align*}
\|F_n(t,v_n(t)) - F_\infty(t)v_n(t)\| & = \|\|u_n\|^{-1}_\infty F(t, \|u_n\|_\infty v_n(t)) - F_\infty(t)v_n(t)\| \\
& \le \|u_n\|^{-1}_\infty (\varepsilon \|u_n\|_\infty \|v_n(t)\|_\alpha + m_\varepsilon) \\
& \le \varepsilon + m_\varepsilon/\|u_n\|^{-1}_\infty \qquad\mathrm{for}\quad n\ge 1, \ t\in[0,+\infty).
\end{align*}
Since $\varepsilon > 0$ is arbitrary, we deduce that
\begin{equation*}
\|F_n(t,v_n(t)) - F_\infty(t)v_n(t)\| \to 0 \qquad\mathrm{as}\quad n\to +\infty, \ \text{uniformly for }t\ge 0.
\end{equation*}
Combining this with \ref{efefeff} implies that
$$v_0(t) = S_A(\lambda_0 t)x_0 + \lambda_0\int_0^tS_A(\lambda_0(t - \tau))F_\infty(\tau)v_0(\tau)\, d\tau \qquad\mathrm{for}\quad t\in[0,+\infty),$$
which is a contradiction because the $T$-periodic problem \ref{nnrr2} does not admit nontrivial solutions. \\[5pt]
\emph{Case B.} Suppose that $\lambda_0 = 0$. Since, for any $n\ge 1$, the map $u_n$ is $T$-periodic
$$x_n= S_A(kT\lambda_n)x_n + \int_0^{kT\lambda_n}S_A(kT\lambda_n - s)F_n(s/\lambda_n,v_n(s/\lambda_n)) \, ds \qquad\mathrm{for}\quad n,k\in\mathbb{N}.$$
For any $n\ge 1$ take $k_n:=[\frac{1}{2\lambda_n}]$. Then $t_n:=k_n T\lambda_n \to T/2$ as $n\to +\infty$.
From Lemma \ref{lem-conv} $(ii)$ and \ref{roroee}, it follows that $(x_n)$ is a relatively compact sequence and therefore, without loss of generality, we can assume that there is $x_0\in X^\alpha$ such that $x_n\to x_0$ as $n\to +\infty$. Hence, by \ref{roronn2}
\begin{align*}
\|v_n(t) - x_n\|_\alpha & = \|S_A(\lambda_n t)x_n - x_n\|_\alpha + \lambda_n\int_0^t\|A^\alpha S_A(\lambda_n(t - \tau))F_n(\tau,v_n(\tau))\| \, d\tau \\
& \le \|S_A(\lambda_n t)x_n - x_n\|_\alpha + \lambda_n^{1 - \alpha}\int_0^t K M_\alpha (t - \tau)^{-\alpha} \,d \tau \\
& \le \|S_A(\lambda_n t)x_n - x_n\|_\alpha
+ \frac{K M_\alpha}{1 - \alpha}\lambda_n^{1 - \alpha}T^{1 - \alpha} \qquad\mathrm{for}\quad t\in[0,T],
\end{align*}
where $K$ is the constant from \ref{roroee} and $M_\alpha$ is the constant from \ref{esti1}. Therefore
\begin{equation}\label{uni}
v_n(t) \to x_0 \qquad\mathrm{as}\quad n\to +\infty, \text{ uniformly for }t\in[0,T].
\end{equation}
Furthermore $\|x_0\|_\alpha = \lim_{n\to +\infty}\|v_n\|_\infty = 1$. Let us write the Duhamel formula \ref{roronn2} in the following form
\begin{align*}
\frac{x_n - S_A(\lambda_n T)x_n}{\lambda_n T} & = \frac{1}{T} \int_0^{T}S_A(\lambda_n (T - \tau)) (F_n(\tau,v_n(\tau)) - F_\infty(\tau)v_n(\tau)) \, d\tau \\
& \qquad + \frac{1}{T} \int_0^{T}S_A(\lambda_n (T - \tau))F_\infty(\tau)v_n(\tau)\,d\tau.
\end{align*}
Proceeding as in the Case A, we can prove that
$$\|F_n(t,v_n(t)) - F_\infty(t)v_n(t)\| \to 0 \qquad\mathrm{as}\quad n\to +\infty, \ \text{uniformly for }t\ge 0.$$
Combining this with \ref{uni}, we infer that
\begin{align*}
& \frac{1}{T} \int_0^{T}S_A(\lambda_n (T - \tau)) (F_n(\tau,v_n(\tau)) - F_\infty(\tau)v_n(\tau)) \, d\tau \to 0 &&  \text{and} \\
& \frac{1}{T} \int_0^{T}S_A(\lambda_n (T - \tau))F_\infty(\tau)v_n(\tau)\,d\tau \to \frac{1}{T} \int_0^{T} F_\infty(\tau)x_0 \,d\tau && \text{as} \ n\to+\infty.
\end{align*}
Therefore, similar reasoning as in the Step 2 of the proof of Theorem \ref{th-cont} leads to $- Ax_0 + \widehat{F}_\infty(x_0) = 0$ which contradicts the assumption because $\|x_0\|_\alpha = 1$. Hence the claim is proved. \\[5pt]
\noindent\emph{Step 2.}  Let us define the homotopy $H:[0,T]\times X^\alpha \times [0,1] \to X$ by
\begin{equation*}
H(t,x,\lambda) :=
\left\{\begin{aligned}
& \lambda F(t,\lambda^{-1} x) && \text{for } t\in[0,T], \ x\in X^\alpha, \ \lambda\in[0,1], \\
& F_\infty(t)x && \text{for } t\in[0,T], \ x\in X^\alpha, \ \lambda = 0. \\
\end{aligned}\right.
\end{equation*}
It is not difficult to check that $H$ is a continuous map. We claim that there is $R_1 > R_0$ such that
\begin{equation}\label{aaaaa3}
-A x + \widehat{H}(x,\lambda) \neq 0 \qquad\mathrm{for}\quad \lambda\in[0,1] \text{ and }x\in D(A) \text{ with }\|x\|_\alpha \ge R_1.
\end{equation}
Otherwise there are sequences $(\lambda_n)$ in $[0,1]$ and $(x_n)$ in $X^\alpha$ such that $-A x_n + \widehat{H}(x_n,\lambda_n) = 0$ and $\|x_n\|_\alpha \to +\infty$ as $n\to+\infty$. Let us denote $z_n:= x_n/\|x_n\|_\alpha$. If $\lambda_n = 0$ for some $n \ge 1$, then $z_n = A^{-1}\widehat{F}_\infty z_n$ which is a contradiction with the assumption of theorem. If $(\lambda_n)$ is contained in $(0,1]$, then
\begin{equation}\label{aaaaa}
z_n = A^{-1}(\lambda_n^{-1}\|x_n\|_\alpha)^{-1}\widehat{F}(\lambda_n^{-1}\|x_n\|_\alpha z_n) \qquad\mathrm{for}\quad n\ge 1.
\end{equation}
Let us note that assumption \ref{nnrr} gives
\begin{equation}\label{aaaaa2}
\lim_{\|x\|_\alpha\to +\infty} \frac{\|\widehat{F}(x) - \widehat{F}_\infty x\|}{\|x\|_\alpha} = 0.
\end{equation}
On the other hand, if we denote $\rho_n:= \lambda_n^{-1}\|x_n\|_\alpha$, then $\rho_n\to +\infty$ as $n\to +\infty$. Therefore, using \ref{aaaaa} and \ref{aaaaa2}, we infer that
\begin{equation}\label{nnnaa}
\begin{aligned}
& \|A^\alpha z_n - A^{-1 + \alpha}\widehat{F}_\infty z_n\| =
\|\rho_n^{-1}A^{- 1 + \alpha}\widehat{F}(\rho_n z_n) - A^{-1 + \alpha}\widehat{F}_\infty z_n\| \\
& \qquad \le \|A^{-1 +\alpha}\| \|\widehat{F}(\rho_n z_n) - \widehat{F}_\infty(\rho_n z_n)\|/\rho_n \to 0 \qquad\mathrm{as}\quad n\to+\infty.
\end{aligned}
\end{equation}
Let us observe that, by \cite[Theorem 1.4.8]{Henry}, the operator $A^{-1 + \alpha}:X\to X$ is compact. Combining this fact with \ref{nnnaa}, implies that the sequence $(z_n)$ is relatively compact in $X^\alpha$. Without loss of generality we can assume that $z_n \to z_0$ for some $z_0\in X^\alpha$ with $\|z_0\|_\alpha = 1$. Since
$$\|z_n - A^{-1}\widehat{F}_\infty z_n\|_\alpha = \|A^\alpha z_n - A^{-1 + \alpha}\widehat{F}_\infty z_n\| \qquad\mathrm{for}\quad n\ge 1,$$ using \ref{nnnaa} again, we have
$z_0 = A^{-1}\widehat{F}_\infty (z_0)$, contrary to assumption of the theorem. Therefore \ref{aaaaa3} follows. \\[5pt]
\noindent\emph{Step 3.} Write $B(0,R_1):= \{x\in X^\alpha \ | \ \|x\|_\alpha < R_1\}$. By Step 2, homotopy invariance of topological degree (see Theorem \ref{th-degal}) and \ref{deg-def1} (with $\mu = 0$) we infer that
\begin{align*}
\mathrm{deg_\alpha}(-A + \widehat F, B(0,R_1)) & = \mathrm{deg_\alpha}(-A x + \widehat{H}(\,\cdot\,,1), B(0,R_1)) \\
& = \mathrm{deg_\alpha}(-A x + \widehat{H}(\,\cdot\,,0), B(0,R_1)) \\
& = \mathrm{deg_\alpha}(-A + \widehat F_\infty, B(0,R_1)) \\
& = \mathrm{deg_{LS}}(I - A^{-1}\widehat F_\infty, B(0,R_1)) = \pm 1,
\end{align*}
where the last equality follows from the fact that the linear operator $A^{-1}\widehat F_\infty:X^\alpha\to X^\alpha$ is compact and has trivial kernel. By Theorem \ref{th-cont}, there is $\lambda_0\in (0,1)$ such that, if $\lambda\in(0,\lambda_0]$, then $\Phi^\lambda_T(x)\neq x$ for $x\in \partial U$ and
\begin{equation}\label{nnn111}
\mathrm{deg_{LS}}(I - \Phi^\lambda_T, U) = \mathrm{deg_\alpha}(-A + \widehat{F}, U) = \pm 1.
\end{equation}
By Step 1 and homotopy invariance of Leray-Schauder degree, one has
$$\mathrm{deg_{LS}}(I - \Phi^1_T, U) = \mathrm{deg_{LS}}(I - \Phi^{\lambda_0}_T, U).$$
Combining this with \ref{nnn111} we infer that $\mathrm{deg_{LS}}(I - \Phi^1_T, U) = \pm 1$, which implies that the equation \ref{nowe} admits a $T$-periodic solution as claimed.
\end{proof}

\section{Applications}

In this section we provide applications of the obtained abstract results to particular partial differential equations. We will assume that $\Omega\subset\mathbb{R}^n$, $n\ge 1$, is an open bounded set with the boundary $\partial\Omega$ of class $C^1$. Let us consider the following partial differential equation
\begin{equation}\label{zag11}
u_t = -\mathcal{A} \, u + f(t, x, u, \nabla u), \qquad t > 0, \ x\in\Omega,
\end{equation}
where $\mathcal{A}$ is a differential operator of the following form
$$\mathcal{A} \bar u (x) = -\sum_{i,j=1}^n D_j(a_{ij}(x)D_i \bar u(x)) \qquad\mathrm{for}\quad \bar u\in C^1(\overline{\Omega}),$$
where the coefficients $a_{ij}$ are of class $C^1(\overline\Omega)$ and
\begin{equation}\label{ellip33}
\sum_{1\le i,j\le n}a_{ij}(x)\xi^i\xi^j \ge c_0 |\xi|^2 \qquad\mathrm{for}\quad x\in\Omega, \ \xi\in\mathbb{R}^n, \ \ \text{where} \ \ c_0 > 0.
\end{equation}
Furthermore we assume that $f:[0,+\infty)\times\Omega\times\mathbb{R}\times\mathbb{R}^n\to\mathbb{R}$ is a continuous map such that the following conditions hold: \\[5pt]
\noindent\makebox[22pt][l]{$(E1)$} \parbox[t][][t]{118mm}{there is a constant $L > 0$ such that, if $t\in[0,+\infty)$ and $x\in\Omega$, then
    \begin{align*}
    |f(t,x,s_1,y_1) - f(t,x,s_2,y_2)| & \le L(|s_1 - s_2| + |y_1 - y_2|)
    \end{align*}
for $s_1,s_2\in\mathbb{R}$ and $y_1,y_2\in\mathbb{R}^n$;}\\[5pt]
\noindent\makebox[22pt][l]{$(E2)$} \parbox[t][][t]{118mm}{there is a continuous function $m:[0,+\infty)\to [0,+\infty)$ such that
    \begin{align*}
    |f(t,x,s,y)| & \le m(t)(1 + |s|),
    \end{align*}
for $t\in[0,+\infty)$, $x\in\Omega$, $s\in\mathbb{R}$ and $y\in\mathbb{R}^n$;}\\[5pt]
\noindent\makebox[22pt][l]{$(E3)$} \parbox[t][][t]{118mm}{there is $T > 0$ such that for any $s\in\mathbb{R}$, $y\in\Omega$ and $y\in\mathbb{R}^n$ one has $$f(t + T, x, s, y) = f(t, x, s, y) \qquad\mathrm{for}\quad t\in[0,+\infty);$$}\\[5pt]
\noindent\makebox[22pt][l]{$(E4)$} \parbox[t][][t]{118mm}{there is a $T$-periodic continuous function $f_\infty:[0,+\infty)\to \mathbb{R}$ such that $$\lim_{|s|\to+\infty} f(t,x,s,y)/s = f_\infty(t)$$
uniformly for $t\in[0,+\infty)$, $x\in\Omega$ and $y\in\mathbb{R}^n$. } \\[5pt]
Let us introduce the abstract framework for the equation \ref{zag11}. To this end, we denote $X:=L^p(\Omega)$, for $p\ge 2$, and define the operator $A_p: X\supset D(A_p)\to X$ by
\begin{equation*}
\begin{aligned}
D(A_p) := W^{2,p}(\Omega) \cap W^{1,p}_0(\Omega), \quad A_p \bar u := \mathcal{A} \bar u \quad \text{for} \ \ \bar u\in D(A_p).
\end{aligned}
\end{equation*}
Then, it is known (see e.g. \cite{DloCho}, \cite{Pazy}, \cite{MR500580}) that $A_p$ is a positive sectorial operator with compact resolvent. If we denote $X^\alpha := D(A_p^\alpha)$, for some $\alpha\in(1/2,1)$, then \cite[Theorem 1.6.1]{Henry} implies that $X^\alpha\subset W^{1,p}(\Omega)$ and the inclusion is continuous. Hence we can define $F\colon [0,+\infty)\times X^\alpha \to X$ as a map given, for any $\bar u\in X^\alpha$, by
\begin{equation*}
    F(t,\bar u)(x) := f(t,x, \bar u(x), \nabla \bar u(x)) \qquad\mathrm{for}\quad t\in [0,+\infty), \ x\in\Omega.
\end{equation*}
We call $F$ \emph{the Nemitskii operator} associated with $f$. Simple calculations shows that $F$ is well-defined, continuous and furthermore, assumptions $(E1)-(E3)$ imply that conditions $(F1)$, $(F2)$ and $(F3)$ are satisfied. Now, we are ready to write the equation \ref{zag11} in the following abstract form
\begin{equation}\label{eqq1}
\dot u(t) = - A_p u(t) + F (t,u(t)), \qquad t > 0.
\end{equation}
We prove the following theorem.
\begin{theorem}
If assumptions $(E1)-(E4)$ are satisfied and the number $$\widehat f_\infty:=\frac{1}{T}\int_0^T f_\infty(\tau) \, d\tau$$ is such that $\left\{\widehat f_\infty + y i \ | \ y \in\mathbb{R} \right\}\cap \sigma(A_p) = \emptyset$, then the equation \ref{eqq1} admits a $T$-periodic mild solution.
\end{theorem}
\begin{proof}
We intend to verify the assumptions of Theorem \ref{th-cont2}. Given $t\in[0,+\infty)$, let us define $F_\infty(t): X^\alpha \to X$ by $$(F_\infty(t)u)(x) = f_\infty(t) u(x) \quad \text{for a.a.} \ \ x\in\Omega.$$ First we prove that
\begin{equation}\label{oioi3}
\lim_{\|u\|_\alpha\to +\infty} \frac{\|F(t,u) - F_\infty(t)u\|}{\|u\|_\alpha} = 0 \qquad \text{uniformly for } \ t\in[0,+\infty).
\end{equation}
To this end let $(t_n)$ in $[0,+\infty)$ and $(u_n)$ in $X^\alpha$ be sequences such that $\|u_n\|_\alpha \to +\infty$ as $n\to +\infty$. Given $n\ge 1$, define $z_n:=u_n/\|u_n\|_\alpha$. Since, by \cite[Theorem 1.4.8]{Henry}, the inclusion $X^\alpha\subset X$ is compact, there is a subsequence $(z_{n_k})$ together with functions $z_0,g\in L^p(\Omega)$ such that $z_{n_k} \to z_0$ in $L^p(\Omega)$, $z_{n_k}(x) \to z_0(x)$ for a.a. $x\in\Omega$ and $|z_{n_k}(x)| \le g(x)$ for $k \ge 1$. If we denote
$$C_k(x):= \mu_{n_k}^{-1}f(t_{n_k}, x, \mu_{n_k}z_{n_k}(x), \mu_{n_k}\nabla z_{n_k}(x)) - f_\infty(t_{n_k}) z_{n_k}(x)$$
where $\mu_n:=\|u_n\|_\alpha$, then, from assumptions $(E2)$ and $(E4)$, one has
\begin{equation}\label{oioi1}
|C_k(x)|^p \to 0  \ \ \text{ for a.e. } x\in\Omega, \ \text{ as } k\to+\infty.
\end{equation}
Furthermore, by assumption $(E2)$, we infer that
\begin{equation*}
|C_k(x)|^p \le (K_1(K_2 + g) + K_3 g)^p \text{ a.e. on } \ x\in\Omega,
\end{equation*}
where $K_1:=\sup_{t\in[0,T]} m(t)$, $K_2:=\sup_{n\ge 1}\mu_n^{-1}$ and $K_3:=\sup_{t\in[0,T]} |f_\infty(t)|$. Combining this inequality with \ref{oioi1} gives \ref{oioi3} as desired. Furthermore, we see that
$$\ker(-A_p + \widehat F_\infty) = \ker(- A_p + \widehat f_\infty I) = \{0\}$$ as $\widehat f_\infty\not\in \sigma(A_p)$. Therefore, it remains to show that if $\lambda\in(0,1]$, then the problem
\begin{equation}\label{nnrr255}
\left\{\begin{aligned}
\dot u(t) & = \lambda(-A_p + F_\infty(t))u(t), \qquad t\in[0,+\infty) \\
u(0) & = u(T)
\end{aligned}\right.
\end{equation}
does not admit nontrivial mild solutions. On the contrary, suppose that there is $\lambda\in(0,1]$ and a nontrivial $T$-periodic mild solution $u$ for the equation \ref{nnrr255}. It is not difficult to check by direct calculations that $$u(t) = \exp\left(\int_0^t\lambda f_\infty(\tau)\, d\tau\right) S_{\lambda A_p}(t)u(0) \qquad\mathrm{for}\quad t\ge 0$$ and therefore, by $T$-periodicity, $$e^{-\lambda T\widehat f_\infty}u(0) = \exp\left(- \lambda \int_0^T f_\infty(\tau)\, d\tau\right) u(0) = S_{\lambda A_p}(T)u(0).$$ Hence,
from \cite[Theorem 16.7.2]{H-F} (see also \cite[Corollary 3.8]{MR1721989}), it follows that
\footnote{\, For $K\subset X$, by $\overline{\mathrm{lin}} \, K$ we denote the closure of the linear space spanned on the set $K$.}
\begin{equation}\label{formop}
u(0)\in \overline{\mathrm{lin}}\bigcup_{k\in\mathbb{Z}} \, \ker\left(\left(\widehat f_\infty + \frac{2k\pi}{\lambda T} i\right) I - A_p\right).
\end{equation}
Since $u$ is nontrivial, we have in particular that $u(0)\neq 0$ and hence, there is $k\in\mathbb{Z}$ such that $\widehat f_\infty + \frac{2k\pi}{\lambda T} i$ is an eigenvalue of $A_p$. But this contradicts the hypotheses and consequently \ref{nnrr255} does not admit nontrivial mild solutions. Therefore the assumptions of Theorem \ref{th-cont2} are satisfied and the equation \ref{zag11} admits a $T$-periodic mild solution as desired.
\end{proof}

\section{Appendix}
In this section we briefly provide a construction and important properties of the topological degree for perturbations of sectorial operators that we exploit in the previous sections. Although the following definition is convenient from the point of view of our studies, the topological degree for perturbations of sectorial operators is actually a special case of the Mawhin's coincidence degree for an abstract operator which is the sum of a Fredholm linear map of index zero and a nonlinear mapping equipped with some compactness properties. For more details on this framework, we refer the reader to \cite{gaines}, \cite{mawh2}.

\begin{definition}
Given a Banach space $X$ and $\alpha\in(0,1)$, {\em the class of admissible operators} is the set $\mathcal{A}(\alpha, X)$ which consists of the maps $-A + F: \overline{U}\cap D(A) \to X$ with the following properties:
\begin{itemize}
\item $A$ is a positive sectorial linear operator and with compact resolvents,
\item the map $F:X^\alpha \to X$ is continuous and bounded, that is, $F(V)$ is a bounded subset in $X$ for every bounded subset $V\subset X^\alpha$,
\item the set $U\subset X^\alpha$ is open bounded and $-A x + F(x)\neq 0$ for $x\in\partial U\cap D(A)$.
\end{itemize}
\end{definition}
Let $-A + F: D(A)\cap \overline{U} \to X $ be a map of class $\mathcal{A}(\alpha,X)$ for some $\alpha\in(0,1)$. Since $A$ is positive, there is constant $\omega > 0$ such that $(-\omega, +\infty)\subset\rho(A)$. Given arbitrary $\mu \in (-\omega, +\infty)\subset\rho(-A)$, we define {\em the topological degree for perturbations of sectorial operators} by
\begin{equation}\label{deg-def1}
\mathrm{deg_\alpha}(-A + F, U) := \mathrm{deg_{LS}}(I - i_\alpha(\mu I + A)^{-1}(\mu I + F), U),
\end{equation}
where $i_\alpha : X^1 \to X^\alpha$ is the continuous inclusion and $\mathrm{deg_{LS}}$ is a Leray-Schauder topological degree. It can be proved that this definition is correct and independent of the choice of the parameter $\mu \in (-\omega, +\infty)$. \\

In the following theorem we collect the expected properties of the map $\mathrm{deg_\alpha}$. For the proof we refer the reader to \cite{gaines}, \cite{mawh2}.

\begin{theorem}\label{th-degal}
The topological degree $\deg_\alpha$ admits the following properties. \\[5pt]
\makebox[8mm][l]{$(D1)$} \parbox[t]{118mm}{(Existence) If $\mathrm{deg_\alpha}(-A + F, U) \neq 0$ then $-A x + F(x) = 0$ for some $x\in U$.} \\[5pt]
\makebox[8mm][l]{$(D2)$} \parbox[t]{118mm}{(Additivity) Let $-A + F: \overline{U}\cap D(A) \to X$ be an admissible map and let $U_1, U_2\subset U$ be disjoint open sets such that $$\{x\in \overline{U}\cap D(A) \ | \ -A x + F(x) = 0\}\subset U_1\cup U_2.$$ Then $\mathrm{deg_\alpha}(-A + F, U) = \mathrm{deg_\alpha}(-A + F, U_1) + \mathrm{deg_\alpha}(-A + F, U_2)$.}\\[5pt]
\makebox[8mm][l]{$(D3)$} \parbox[t]{118mm}{(Homotopy invariance) Let $F:[0,1]\times X^\alpha \to X$, where $\alpha\in(0,1)$, be a continuous map, transforming bounded sets in $[0,1]\times X^\alpha$ onto bounded sets in $X$. If $U\subset X^\alpha$ is an open bounded set such that $-A x + F(\lambda,x)\neq 0$ for $\lambda\in [0,1]$ and $x\in\partial U\cap D(A)$, then
\begin{align*}
\mathrm{deg_\alpha}(-A + F(0,\,\cdot \,), U) = \mathrm{deg_\alpha}(-A + F(1,\,\cdot \, ), U).
\end{align*}}
\makebox[8mm][l]{$(D4)$} \parbox[t]{118mm}{(Normalization) Let $A:X\supset D(A) \to X$ be a positive sectorial operator with compact resolvents. If $x_0 \in A(U\cap D(A))$, then the map $-A + x_0: \overline{U}\cap D(A) \to X$ is admissible and \ $\mathrm{deg_\alpha}(-A + x_0, U) = 1$.}
\end{theorem}


\def\cprime{$'$} \def\polhk#1{\setbox0=\hbox{#1}{\ooalign{\hidewidth
  \lower1.5ex\hbox{`}\hidewidth\crcr\unhbox0}}} \def\cprime{$'$}

\parindent = 0 pt

\end{document}